\documentclass[preprint,12pt]{elsarticle}

\usepackage{amssymb}

\usepackage{amsmath}
\usepackage{bm}

\DeclareMathOperator{\gd}{gd}
\DeclareMathOperator{\sech}{sech}
\DeclareMathOperator{\arsinh}{arsinh}
\DeclareMathOperator{\artanh}{artanh}
\DeclareMathOperator{\kexp}{\exp_{\kappa}}
\DeclareMathOperator{\kln}{\ln_{\kappa}}
\DeclareMathOperator{\ku}{ u_{\kappa} }
\DeclareMathOperator{\Christ}{\Gamma}
\newcommand{\dbar}{d\hspace*{-0.08em}\bar{}\hspace*{0.1em}}

\journal{Physica A}

\begin{document}

\begin{frontmatter}



\title{On some information geometric structures concerning Mercator projections}


\author{Tatsuaki Wada}

\address{Region of Electrical and Electronic Systems Engineering, Ibaraki University, 316-8511, Japan}

\begin{abstract}
Some information geometric structures concerning the Mercator projections are studied.
It is known that a loxodrome on the surface of the globe is related to the straight line 
on a Mercator map by the Mercator projection. It is not well known that an affine connection with torsion plays a fundamental role to describe an auto-parallel path on the surface. Based on these information geometric structures, Gauss distribution is reconsidered from the view point of the affine connection with a torsion.
Some relations with deformed functions are also pointed out.
\end{abstract}

\begin{keyword}
information geometry \sep affine connections \sep torsions \sep deformed exponentials



\end{keyword}

\end{frontmatter}



\section{Introduction}
\label{intro}
The method of information geometry \cite{Amari} applies some techniques of affine differential geometry to various fields related with probability distributions, and
provides a useful tool to study several applications in statistics, machine learning, and statistical physics.  
A distinct characteristics of the geometrical structures of information geometry is flatness, which is characterized by a zero affine curvature or affine geodesic.
With an appropriate affine connection $\nabla$, a generalization of the straight line on an affine manifold is determined as an auto-parallel curve, which is called \textit{affine geodesic}.  
Indeed, for an affine parametrized curve $\gamma(t)$ on a smooth manifold $\mathcal{M}$ with an affine connection $\nabla$, the geodesic equation is given by
\begin{equation}
  \nabla_{\frac{d \gamma}{dt}} \frac{d \gamma(t)}{d t} = 0,
\end{equation}
whose solution is auto-parallel, i.e., its tangent vector at a point $t$ of the curve $\gamma(t)$ is parallel transported with respect to itself when it is transported along  $\gamma(t)$.

On the one hand, some recent developments in information geometry are affected \cite{MAGE} by the progress of the generalized thermostatistics \cite{Naudts} based on the generalized entropies  and the deformed exponential functions such as $q$-deformed \cite{Tsallis}, $\kappa$-deformed \cite{kappa}, and Newton's deformed \cite{Newton} exponential functions. 
In these developments one of the remarkable points is \textit{conformal flattening}, in which a flat structure is obtained by applying conformal transformation with an appropriate conformal factor \cite{OMA10,OMA12,WM17}. For the details of the conformal flattening on the probability simplex, please consult Ref. \cite{O18}

On the other hand,
the Mercator projection is a cylindrical projection and most popular for us through the standard maps in the world atlas.
\begin{figure}[h]
\begin{center}
\includegraphics{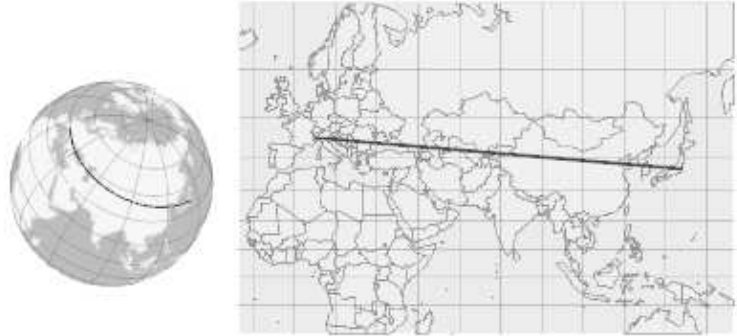}
\end{center}
\caption{Globe of the earth and a Mercator map. A straight line on the Mercator map is the Marcator projection of a loxodrome on the surface of the globe.}
\label{fig:Mercator}       
\end{figure}
Though Mercator projection heavily distorts the size and shape of objects as distance from the equator increases, at every point the scale factor of east-west stretching is same as that of north-south stretching. Consequently it is a conformal projection.
Historically it became the standard map projection for nautical purposes because of its ability to represent lines of constant course, known as  \textit{rhumb lines} or \textit{loxodromes}, as straight segments that conserve the angles $\phi$ with the meridians.
As shown in Fig. \ref{fig:Mercator}, the Mercator projection of a loxodrome on the sruface of the globe is a straight line on a Mercator map. This means that a loxodrome is an affine geodesic (auto-parallel path) 
for an appropriate affine connection.
In this way, the two properties, \textit{conformal transformation} and \textit{auto-parallel path} are common important geometrical ingredients in both information geometry and Mercator projections.

The purpose of this contribution is to explore the information geometric structures concerning Mercator projections. Based on these information geometric structures, Gauss distribution, which
is a fundamental model in information geometry, is reconsidered and it can be characterized by utilizing
an affine connection with torsion.
The rest of the paper is organized as follows.
The next section briefly review the Mercator projection, in which Gudermannian function is naturally arises in the inverse equation of the Mercator projection. Section 3 considers some information geometric structures concerning the Mercator projections. Loxodrome on a sphere is anto-parallel path with respect to an affine connection with torsion. In subsection 3.1 it is shown that a loxodorome
on a pseudosphere is considered as an auto-parallel path which characterizes
  Gauss distribution. Section 4 describes some papers which discussed torsions in the field of information geometry, and point out some relations with deformed functions.
 Final section is devoted to conclusions and perspectives.
Appendix A provides a short review of anholonomic coordinates transformation and Appendix B is
a short explanation of the metric and coefficients of a connection under the conformal transformation.

\section{Mercator projection and Gudermannian function}
\label{sec:1}

The Gudermannian function \cite{Gfunc} is defined for a real variable $x$ by
\begin{equation}
 \gd(x) := \int_0^x \frac{ds}{\cosh (s)},
  \label{gd}
\end{equation}
which arises in the inverse mapping of the Mercator projection.
The Gudermannian function has different expressions and, for example, it is also expressed as
\begin{equation}
  \gd(x) = \arctan \big( \sinh (x) \big).
\end{equation}
The inverse function of  the Gudermannian function is given by
\begin{equation}
  \gd^{-1} (\theta) = \int_0^{\theta} \frac{ds}{\cos (s)},  \quad -\frac{\pi}{2} < \theta < \frac{\pi}{2}.
\end{equation}
The derivative of $\gd(x)$ and that of $\gd^{-1}( \theta)$ are
\begin{equation}
 \frac{d}{dx} \gd (x) = \sech (x), \quad \frac{d}{d \theta} \gd^{-1} (\theta) = \frac{1}{\cos (\theta)},
\end{equation}
respectively.
It is worth noting that the Gudermannian functions connect the trigonometric and hyperbolic functions, e.g.,
\begin{align}
 \sin \big( \gd(x) \big)	&=	\tanh (x), \quad \cos \big( \gd(x) \big)	=	\sech(x),	\\
 \tan \big( \gd(x) \big) &=	\sinh(x),
\end{align}
without using an imaginary number.

Mercator projection \cite{Mercator} transforms a point $(\varphi, \theta)$ on the surface of the globe to
a point $(x, y)$ on the plane of Mercator map as
\begin{equation}
  x = R \, (\varphi - \varphi_0), \quad  y = R \, \gd^{-1} (\theta),
\label{MercatorMap}
\end{equation}
where $R$ is the radius of the globe, $\theta$ is the latitude ($-\pi /2 < \theta < \pi /2$), $\varphi$ is the longitude ($-\pi \le \varphi \le \pi$), 
and $\varphi_0$ is the longitude of an arbitrary central meridian (usually that
of Greenwich is set $\varphi_0=0$).
Mercator projection is a conformal transformation. Indeed, we see that since the infinitesimal area $dA := R \cos(\theta) d\varphi \wedge R d \theta$ on a surface of the globe is transformed to $dx \wedge dy = R d\varphi \wedge R d\theta / \cos(\theta)$ on the plane of the Mercator map, the conformal factor of this transformation is $1/\cos(\theta)$.

\begin{figure}[h]
\begin{center}
\includegraphics[width=0.4\textwidth]{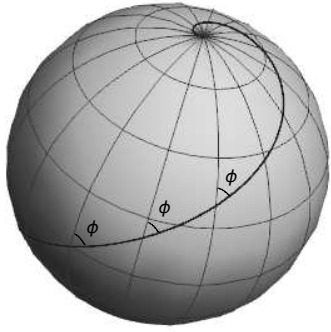}
\end{center}
\caption{A loxodrome on a sphere. It cuts all meridians at a constant angle $\phi$.}
\label{fig:loxodrome}       
\end{figure}

The inverse transformation of the Mercator projection is given by
\begin{equation}
  \varphi = \varphi_0 + \frac{x}{R}, \quad \theta = \gd \left(\frac{y}{R} \right).
\label{invMercator}
\end{equation}
Loxodrome is a curve (or path) which cuts all meridians on a given surface at a constant angle  $\phi$ ($0<\phi<\pi$ but $\phi \ne \pi/2$) as shown in Fig. \ref{fig:loxodrome}.
The equation of a loxodrome on the surface of the globe is written by
\begin{align}
  \theta = \gd \big( \cot(\phi) \, (\varphi-\varphi_0) \big), \quad
\textrm{ or } \quad
 \varphi = \varphi_0 + \tan(\phi) \, \gd^{-1}( \theta). 
\label{loxo}
\end{align}

\section{Information geometric structures concerning the Mercator projections}
\label{sec:3}

Having described the basic of the Mercator projections, we now consider some information geometric structures such as affine geodesics (auto-parallel paths), and conformal flattening with an appropriate conformal factor.
In the following we find an affine connection for which a loxodrome on the sphere with radius $R$ is a auto-parallel path (affine geodesics). Such a path can be considered as an orbital along which a free particle moves. This kind of view is the heart of Einstein's general relativity \cite{Einstein}. We can consider an affine geodesic
as a "straightest" line along which a free particle moves.

Now we begin with a straight path on a Mercator map, and consider its mapping by the inverse Mercator projection. 
Suppose we have a free particle with unit mass, whose dynamical simultaneous equations are
\begin{equation} 
   \frac{d^2}{ dt^2} \, x(t) = 0, \quad \frac{d^2}{ dt^2} \, y(t) =0,
\end{equation}
on a two-dimensional plane ($x\,y$-plane). The solution can be considered as a trajectory along a straight line on a Mercator map.
\begin{equation}
  x(t) =  \tan(\phi) \, t, \quad y(t) =   t,
 \label{line}
\end{equation}
where $\phi$ is a constant angle which determines a magnitude of the velocity $v$ of the particle moving along
this straight line, i.e.,
\begin{align}
  v_x &= \frac{d}{d t} x(t) = \tan(\phi), \quad v_y = \frac{d}{dt} \, y(t) = 1, \textrm{ and }, \nonumber \\
v &= \sqrt{v_x^2 + v_y^2} = \sqrt{\tan^2(\phi) + 1} = \frac{1}{\cos(\phi)}.
\end{align}
Soon later one will see that this $\phi$ is equal to the angle at which the loxodrome cuts each meridian
as shown in Fig. \ref{fig:loxodrome}. 
Recall that this straight line \eqref{line} and a loxodrome \eqref{loxo} are related by the inverse Mercator projection \eqref{invMercator} with the conformal factor of $\cos(\theta)$.
Taking the derivative of \eqref{MercatorMap} and multiplying the conformal factor $\cos(\theta)$, we have
\begin{equation}
    \dbar x = \cos(\theta) R \, d\varphi, \quad dy = R \, d\theta.
\end{equation}
This is an anholonomic coordinate transformation from the geographic coordinate system $\bm{q} := (\varphi, \theta)$ to $\bm{r} := (x, y)$ 
of the following form
\begin{equation}
     \dbar r^i =   e^i{}_{\mu} \, d q^{\mu}, \qquad i=x, y, \quad \mu=\varphi, \theta,
  \label{anholonomic}
\end{equation}
with
\begin{align}
e^i{}_{\mu} =
\begin{pmatrix}
   R \cos (\theta) & 0 \\
   0 & R
\end{pmatrix}
.
\end{align}
Note that since the exterior derivative of $\dbar x$ is non-zero as shown by
\begin{equation}
  d \, (\dbar x) = -R \sin(\theta) \, d\theta \wedge d \varphi = \frac{\tan(\theta)}{R} \, \dbar x \wedge d y,
\end{equation}
$\dbar x$ is an inexact differential. In order to distinguish it from an exact differential $d$, we use
the notation $\dbar$ through out this paper.

The non-zero matrix elements are $e^x{}_{\varphi}= R \cos(\theta)$ and $e^y{}_{\theta}=R$, and
the inverse ${\rm E}_i{}^{\mu}$ of $e^i{}_{\mu}$ is
\begin{equation}
{\rm E}_i{}^{\mu} =
\begin{pmatrix}
   \frac{1}{R \cos (\theta)} & 0 \\
   0 & \frac{1}{R}
\end{pmatrix}.
\end{equation}
Then after straight-forward calculations and using \eqref{Wconnection}, we find that the non-zero coefficient of the affine connection $\nabla$ is
\begin{equation}
{\Christ}^{\varphi}{}_{\theta \varphi} = {\rm E}_x{}^{\varphi} \, \frac{\partial}{\partial \theta} e^x{}_{\varphi} = \frac{1}{R \cos(\theta)} \, \frac{\partial}{\partial \theta}  R \cos(\theta)
= - \tan(\theta).
\end{equation}
Since the other coefficients of the affine connection $\nabla$ are zero, it follows that
\begin{equation}
  {\rm T}^{\varphi}{}_{\theta \varphi} := {\Christ}^{\varphi}{}_{\theta \varphi} - {\Christ}^{\varphi}{}_{\varphi \theta} = \tan(\theta) \ne 0,
\end{equation}
which states that this connection $\nabla$ has a torsion.


\subsection{Gauss  distribution}
We here discuss a Gauss distribution, which is a fundamental statistical model in information geometry \cite{Amari}, from a different point of view. 
The probability density function (pdf) of Gauss  distribution, or normal distribution $N(\mu, \sigma^2)$, is
\begin{equation}
  p(x; \mu, \sigma) := \frac{1}{\sqrt{2 \pi \sigma^2}} \, \exp \left[ -\frac{(x - \mu)^2}{2 \sigma^2} \right],
\end{equation}
where $x$ is a random variable, $\mu$ and $\sigma^2$ stand for the mean and variance, respectively.
The set of all  Gauss distributions is a two-dimensional manifold with the coordinate system $(\mu, \sigma)$ but $\sigma > 0$.
It is well known \cite{Amari} that Fisher-Rao metric of this pdf is
\begin{equation}
  g^{\rm FR}_{ij}(\mu, \sigma) =
 \frac{1}{\sigma^2}
\begin{pmatrix}
   1 & 0 \\
   0 & 2
\end{pmatrix}  
.
\end{equation}
As a result the line element $ds$ satisfies
\begin{equation}
  ds^2  = \frac{ d \mu^2 + 2 d \sigma^2}{\sigma^2},
\end{equation}
which is the metric of Poincar\'e upper half-plane $\mathbb{H}=\{(\mu, \sigma) \vert \sigma>0 \}$ with the negative curvature of $-1/2$.

A pseudosphere with a pseudoradius $R$ has a constant curvature $-1/R^2$ and it is known
that the portion $\sigma > 1$ and $-\pi \le \mu \le \pi$ of the upper half-plane $\mathbb{H}$ is
mapped to the surface of a pseudosphere. As a result, only the limited area in the upper half-plane
is mapped to the surface of a pseudosphere.
Practically we can overcome this difficulty by making the suitable normalization $\tilde{\sigma} := \sigma / \sigma_{\rm min}$ and $\tilde{\mu} := \pi \mu / \vert \mu_{\rm max} \vert$
with an appropriate minimum value $\sigma_{\rm min}$ and an appropriate maximum absolute value $\vert \mu_{\rm max} \vert$. In this way, we use an appropriate normalized set $(\tilde{\mu}, \tilde{\sigma})$ hereafter, i.e.,
\begin{equation}
  ds^2  = \frac{ d \tilde{\mu}^2 + 2 d \tilde{\sigma}^2}{\tilde{\sigma}^2}.
 \label{Poincare}
\end{equation}

There are different expressions for describing a pseudosphere, and a Cartesian parametrization
of them is
\begin{align}
  \big( & r(u)  \cos(\varphi), r(u)  \sin(\varphi),  h(u) \big), \quad
\textrm{ with }
r(u) = R \exp\left (-\frac{u}{R} \right), \textrm{ and }\nonumber \\
& h(u) = R \artanh \left( \sqrt{1-\exp\left(-2\frac{u}{R} \right) } \right)-R\sqrt{1-\exp\left(-2\frac{u}{R} \right)},
\end{align}
where $0 \le u<\infty$ and $0 \le \varphi \le 2 \pi$.
The corresponding Riemann metric is
\begin{equation}
  ds^2 = R^2 \exp \left( -\frac{2u}{R} \right) d \varphi^2 +  d u^2.
 \label{pSphere}
\end{equation}
Recalling $R=\sqrt{2}$, and by comparing \eqref{Poincare} and \eqref{pSphere}, we find that
$du = \sqrt{2} d \tilde{\sigma} / \tilde{\sigma}$ and $\sqrt{2} \exp(-u/\sqrt{2}) d \varphi = \sqrt{2} d \varphi / \tilde{\sigma} = d \tilde{\mu} / \tilde{\sigma}$. As a result we can relate them as
\begin{equation}
  \varphi = \frac{\tilde{\mu}}{R}= \frac{\tilde{\mu}}{\sqrt{2}}, \quad 
 u = R \ln \tilde{\sigma} = \sqrt{2} \, \ln \tilde{\sigma}.
 \label{musigma}
\end{equation}

As similar as in the case of the sphere, we can consider the following coordinate transformations.
\begin{equation}
    \dbar x = R \exp\left(-\frac{u}{R} \right)  \, d\varphi, \quad dy = du.
  \label{AHtrans}
\end{equation}
This is an anholonomic coordinate transformation from the coordinate system $\bm{q} := (\varphi, u)$ to $\bm{r} := (x, y)$ 
of the following form
\begin{equation}
     \dbar r^i =   e^i{}_{\nu} \, d q^{\nu}, \qquad i=x, y, \quad \nu=\varphi, u,
\end{equation}
with the zweibein
\begin{equation}
e^i{}_{\nu} =
\begin{pmatrix}
   R \exp\left(-\frac{u}{R} \right) & 0 \\[1ex]
   0 & 1
\end{pmatrix}  
.
\end{equation}
Since the exterior derivative of $\dbar x$ is non-zero as shown by
\begin{align}
  d \, (\dbar x) &= -\exp\left(-\frac{u}{R} \right) \, du \wedge d \varphi = 
-\frac{1}{R} \, du \wedge \left(R \exp\left(-\frac{u}{R} \right)  d \varphi \right) 
\nonumber \\
&= \frac{1}{R} \, \dbar x \wedge d y,
\end{align}
$\dbar x$ is an inexact differential. 

The non-zero matrix elements are $e^x{}_{\varphi}= R \exp(-u/R)$ and $e^y{}_{\theta}=1$, and
the inverse ${\rm E}_i{}^{\nu}$ of $e^i{}_{\nu}$ is
\begin{equation}
{\rm E}_i{}^{\nu} =
\begin{pmatrix}
   \frac{1 }{R}   \exp\left(\frac{u}{R} \right) & 0 \\[1ex]
   0 &  1
\end{pmatrix}  
.
\end{equation}
Then after straight-forward calculations and using \eqref{Wconnection}, we find that the non-zero coefficient of the affine connection $\nabla^{\rm ps}$ on the surface of pseudosphere is
\begin{equation}
{\Christ}^{\varphi}{}_{u \varphi} = {\rm E}_x{}^{\varphi} \, \frac{\partial}{\partial u} e^x{}_{\varphi} =  \frac{1 }{R}   \exp\left(\frac{u}{R} \right) \frac{\partial}{\partial u}  R  \exp\left(-\frac{u}{R} \right)
= - \frac{1}{R}.
\end{equation}

Next we confirm a loxodrome on the pseudosphere is an auto-parallel path with respect to the affine connection $\nabla^{\rm ps}$. 
Let us setting
\begin{equation}
  u(t) = t,
\end{equation}
then a loxodorome on the pseudosphere is written by
\begin{equation}
  \varphi(t) = \varphi_0 + \tan(\phi) \, \exp\left(\frac{t}{R} \right)
\end{equation}
Taking derivative with respect to $t$, we have
\begin{equation}
  \frac{d \varphi(t)}{d t}  = \frac{\tan(\phi)}{R} \, \exp\left(\frac{t}{R} \right), \quad
  \frac{d^2 \varphi(t)}{d t^2}   = \frac{\tan(\phi)}{R^2} \, \exp\left(\frac{t}{R} \right).
\end{equation}
Consequently we confirm that the auto-parallel equation is satisfied as follows. 
\begin{align}
  \frac{d^2 \varphi(t) }{dt^2} 
   &+ \Christ^{\varphi}{}_{\varphi u} \, \frac{d \varphi(t)}{dt} \, \frac{d u(t) }{dt}  
\nonumber \\
 &= \frac{\tan(\phi)}{R^2} \, \exp\left(\frac{t}{R} \right) - \frac{1}{R} \,  \frac{\tan(\phi)}{R} \, \exp\left(\frac{t}{R} \right) =
   0.
\end{align}
\begin{figure}[h]
\begin{center}
\includegraphics[width=0.4\textwidth]{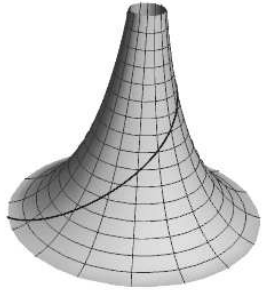}
\end{center}
\caption{A loxodrome on a pseudosphere. }
\label{fig:loxodrome_on}       
\end{figure}

Next by making a conformal transformation of \eqref{AHtrans} by multiplying the conformal factor of
$\exp(u/R)$, we have
\begin{equation}
    d \tilde{x} = R  \, d\varphi, \quad d \tilde{y} = \exp\left(\frac{u}{R} \right)  du,
  \label{conformal}
\end{equation}
which describe Euclidean plane $(\tilde{x}, \tilde{y})$.
Indeed, by using \eqref{UnderCT} with $\lambda = u / R$, we confirm that the conformal
transformed metric and coefficient of the affine connection are
\begin{align}
 \tilde {g}_{ij} &= \exp\left( 2 \frac{u}{R} \right) \, g_{ij} =
\begin{pmatrix}
   R & 0 \\
   0 &  \exp\left(\frac{u}{R} \right)
\end{pmatrix}  
, \\
  \tilde{\Christ}\,{}^{\varphi}{}_{u \varphi} &= \Christ{}^{\varphi}{}_{u \varphi} 
+ \frac{\partial}{\partial u} \, \frac{u}{R} 
   = -\frac{1}{R} + \frac{1}{R} = 0, \\
 \tilde{\Christ}\,{}^{\varphi}{}_{\varphi u} &= \Christ{}^{\varphi}{}_{\varphi u} 
+ \frac{\partial}{\partial u} \, \frac{u}{R} 
   = 0 + \frac{1}{R} = \frac{1}{R},
\end{align}
respectively.
With this conformal mapping, the loxodorome on the pseudosphere is mapped to a straight line
in the $\tilde{x} \, \tilde{y}$-plane.
Substituting \eqref{musigma} into \eqref{conformal}, we finally obtain
\begin{equation}
    d \tilde{x} = R  \, \frac{d \tilde{\mu}}{R} = d\tilde{\mu}, \quad 
    d \tilde{y} = \tilde{\sigma}  d (R \ln \tilde{\sigma}) = R d \tilde{\sigma}
= \sqrt{2} d \tilde{\sigma}.
\end{equation}
Hence this Euclidean plane $(\tilde{x} = \tilde{\mu}, \tilde{y} = \sqrt{2} \tilde{\sigma})$ is spanned by the normalized- mean $\tilde{\mu}$ and standard deviation $\tilde{\sigma}$ of Gauss pdf, and a straight line in this plane characterizes the statistical model of Gauss pdf.

\section{Discussions}
\label{dis}

As well known, in the standard method of information geometry \cite{Amari}, torsion-free
connections are assumed from the onset of the theory.  Consequently curvatures play a dominant role,
as in the case of standard general relativity \cite{Einstein}, in which Levi-Civita connection is used and it is the torsion-free metric connection, i.e., the torsion-free connection on the tangent bundle preserving a given Riemannian metric.
In contrast, the situation is different in the case discussed in section \ref{sec:3}. An affine connection  with torsion plays an important role, which reminds us a teleparallel \cite{teleparallel} approach to general relativity.
In the literature on information geometry, there exist a few papers which discuss non-zero torsion.
To the best of my knowledge, it was Kurose \cite{Kurose07} that considered torsions in statistical manifold for the first time. Matsuzoe \cite{M07,M10} developed the statistical manifold admitting torsion (SMAT),
and Henmi \cite{Henmi17} applied it to the field of statistics. It is hence very interesting to further study SMAT from the view point of the results obtained in this paper for future research works.

Next I'd like to point out some relations with a deformed function. For this purpose, let us briefly review
the $\kappa$-deformed functions \cite{kappa}. 
For a positive real parameter $\kappa$, the deformed exponential function defined by
\begin{equation}
   \kexp( x) := \left[ \kappa x + \sqrt{1 + \kappa^2 x^2} \right]^{\frac{1}{\kappa}},
\end{equation}
is called $\kappa$-exponential function. The inverse function is called $\kappa$-logarithmic function, which is written by
\begin{equation}
   \kln( x) := \frac{x^\kappa -  x^{-\kappa}}{2 \kappa}.
  \label{klog}
\end{equation}
They are fundamental ingredients in one  \cite{WMS15} of the recent developments of information geometry.

Now we focus on the Mercator mapping \eqref{MercatorMap}.
Introducing a new parameter $\chi = R \tan \theta$, it follows that
\begin{equation}
  y(\chi) = R \gd^{-1}(\theta) = R \arsinh \big( \tan(\theta) \big) 
  = R \arsinh \left( \frac{\chi}{R} \right).
\label{k-map}
\end{equation}
Setting $\kappa = 1/R$ we find that
\begin{equation}
  y(\chi) = \frac{1}{\kappa} \arsinh ( \kappa \chi ) 
    = \frac{1}{\kappa} \ln \left( \kappa \chi + \sqrt{1 + \kappa^2 \chi^2 } \right)
   = \ln \left( \kexp (\chi) \right),
\end{equation}
which reduces to the identical mapping $y=\chi$ in the limit of $\kappa \to 0$.
Hence a variant \eqref{k-map} of the Mercator mapping is considered as a $\kappa$-deformation.
In addition, the geometrical meaning of the deformed parameter $\kappa$ in this mapping is very clear as $\kappa=1/R$.
The standard undeformed case of  $\kappa=0$ is corresponding to the limit of the radius $R \to \infty$. 

Such a deformed function is useful for making an immersion (or representation) \cite{WM17} to a statistical manifold in information geometry. Indeed Naudts \cite{Naudts} introduced the generalization of the standard logarithmic function, which is called
the $\phi$-logarithmic function. It is defined by
\begin{equation}
  \ln_{\phi} (x) := \int_1^x \frac{ds}{\phi(s)}, 
\end{equation}
for a positive increasing function $\phi(s)$. For example, choosing 
\begin{equation}
 \phi(s) = \frac{2s}{s^{\kappa}+s^{-\kappa}},
\end{equation}
leads to the $\kappa$-logarithmic function \eqref{klog}.

Now let us introduce another $\kappa$-deformed function defined by
\begin{equation}
 \ku(x) := \frac{x^{\kappa} + x^{-\kappa}}{2}
 = \cosh \big[ \kappa \ln(x) \big],
\end{equation}
which reduces to unit function in the limit of $\kappa \to 0$.
Then since $x \, \ku(x) > 0$ and
\begin{equation}
 \frac{d}{d x} x \, \ku(x) = \cosh \big[ \kappa \ln(x) \big] + \kappa \sinh\big[ \kappa \ln(x) \big] 
= \ku(x) + \kappa^2 \kln(x) > 0,
\end{equation}
for a positive $x$, we see that $x \, \ku(x)$ is a positive increasing function for a positive $x$.
Then we consider the following Naudts $\phi$-logarithmic function with this positive increasing function $\phi(x) = x \, \ku(x)$, i.e.,
\begin{equation}
  \ln_{\phi} (x) := \int_1^x \frac{ds}{s \, \ku(s)}.
  \label{philn}
\end{equation}
Changing the variable $z = \kappa \ln(s)$, it follows that
\begin{equation}
 \ln_{\phi} (x) = \frac{1}{\kappa} \int_1^{\kappa \ln (x)} \frac{dz}{\cosh (z)}
 =  \frac{1}{\kappa} \gd \big( \kappa \ln (x) \big),
  \label{phi-log}
\end{equation}
where we used \eqref{gd} of the Gudermannian function.
Recall the inverse transformation \eqref{invMercator} of the Mercator projection,
\begin{equation}
   R \, \theta = R \gd \left( \frac{y}{R} \right),
  \label{invM}
\end{equation} 
and $\kappa = 1/R$, we can consider the inverse of Mercator projection as a
deformation of the standard logarithmic function a la Naudts. 
This correspondence between \eqref{phi-log} and \eqref{invM} provides the geometric meaning
of the deformed logarithmic function \eqref{philn}.

\section{Conclusions and perspectives}
Some information geometric structures concerning the Mercator projections were discussed.
The Mercator projections are conformal and as pointed out in Introduction, a remarkable key point in recent developments of information geometry is \textit{conformal flattening}, i.e., a geometrical method obtaining a flat structure by applying conformal transformation with an appropriate conformal factor. 
Recently Agricola and Their \cite{AT04} discussed the geodesics of metric connections with vectorial torsion and generalize the Mercator projection. Their results are interesting for further developing the information geometry from the different point of view by using an affine connection with torsion. An attempt for such development has done in section 3.1 for Gauss pdf.

A variant of the Mercator projection is related with the $\kappa$-deformation, where the deformed parameter $\kappa$ has a clear geometrical meaning of $\kappa = 1/R$.
In addition it is related with a $\phi$-logarithmic function by Naudts \cite{Naudts}.

It is known that a loxodrome on the surface of the globe and the straight line on a Mercator map is related by the Mercator projection. However, it is not well known that an affine connection with torsion plays a fundamental role to describe such a loxodrome as an auto-parallel path on the sphere. 
Based on this fact, it is worth to generalize the method of information geometry by incorporating
an torsion such as the works\cite{Kurose07,M07,M10,Henmi17} on SMAT.
Further studies are needed in order to develop a generalization of the information geometry by using an affine connection with torsion.

\section*{Acknowledgement}

The author thanks to Antonio M. Scarfone for useful discussions.
This work is partially supported by Japan Society for the Promotion of Science (JSPS) Grants-in-Aid for Scientific Research (KAKENHI) Grant Number JP17K05341.
The author acknowledges the anonymous referee for his valuable remark.

\appendix


\section{Anholonomic coordinates transformation }
\label{appendix:A}

Consider a free particle in $n$-dimensional Euclidean space, whose dynamic equation is
\begin{equation}
  \frac{d^2 r^i(t)}{dt^2} = 0, \quad i=1,2, \ldots, n.
  \label{Fp}
\end{equation}
Performing the transformation
\begin{equation}
  \dbar r^i =  e^i{}_{\mu} \big ( \bm{q}(t) \big) \, d q^{\mu}(t),
  \label{anholonomicT}
\end{equation}
where the new coordinate $\bm{q} = \{q^{\mu} \}, \mu=1,2, \ldots, n$ and 
\begin{equation}
 e^i{}_{\mu} \big ( \bm{q}(t) \big) := \frac{ \partial r^i}{\partial q^{\mu} },
\end{equation}
are called \textit{vielbein}, and assumed to be the elements of a nonsingular matrix with non-zero determinant.  
If the exterior derivative of the derivative $d r^i$ is non-zero, i.e,
\begin{equation}
   d \, \dbar r^i = d \left( e^i{}_{\mu}(\bm{q}) \, dq^{\mu} \right) \neq 0,
 \label{AHcond}
\end{equation}
the transformation \eqref{anholonomicT} is called \textit{anholonomic} (or non-holonomic) and the $\bm{r} = \{ r^i \}$ is an anholonomic coordinate, and in this case we denote $\dbar r^i$.
 For a holonomic coordinate, say $\{ q^{\mu} \}$, the exterior derivative of $dq^i$ is of course zero.

Now substituting \eqref{anholonomicT} into \eqref{Fp} gives
\begin{equation}
  \frac{d^2 r^i(t)}{dt^2} =  e^i{}_{\mu} ( \bm{q} ) \, \frac{d^2 q^{\mu}(t)}{dt^2} 
  + \frac{\partial}{\partial q^{\mu}} e^i{}_{\mu} ( \bm{q} )   \, \frac{d q^{\nu}}{dt}  \frac{d q^{\mu}}{dt} = 0.
\end{equation}
Multiplying the inverse matrix ${\rm E}^{\rho}{}_i( \bm{q} )$ of $ e^i{}_{\rho} ( \bm{q} ) $ it follows
\begin{equation}
  \frac{d^2 q^{\rho}}{dt^2} 
  + \Christ^{\rho}{}_{\nu \mu} \, \frac{d q^{\nu}}{dt} \frac{d q^{\mu}}{dt} = 0.
\end{equation}
Since this is an auto-parallel equation, the trajectories of the particle in the space $\bm{q}$ are auto-parallels. Here
\begin{equation}
 \Christ^{\rho}{}_{\mu \nu} :=  {\rm E}^{\rho}{}_i \, \frac{\partial}{\partial q^{\mu}} e^i{}_{\nu},
 \label{Wconnection}
\end{equation}
is the coefficient of Weizenb\"ock connection $\nabla^{(W)}$, and they satisfy
\begin{eqnarray}
  \nabla^{(W)}_{\mu} \, {\rm E}_{\nu}{}^i &= \frac{\partial}{\partial q^{\mu}} \, {\rm E}_{\nu}{}^i + \Christ^{\rho}{}_{\mu \nu} \, {\rm E}_{\rho}{}^i =0, \\
\nabla^{(W)}_{\mu} \, {\rm e}^i{}_{\nu} &= \frac{\partial}{\partial q^{\mu}} \, {\rm e}^i{}_{\nu} - \Christ^{i}{}_{\mu \nu} \, {\rm e}^{\rho}{}_i =0,
\end{eqnarray}
respectively.
The Weizenb\"ock connection $\nabla^{(W)}$ is the affine flat connection with zero Cartan curvature and nonzero torsion. The torsion tensor
\begin{equation}
{\rm T}^{\lambda}{}_{\mu \nu}  
  := \Christ^{\rho}{}_{\mu \nu} - \Christ^{\rho}{}_{\nu \mu} \neq 0,
\end{equation}
because of the anholonomic condition \eqref{AHcond}, and
the Riemann curvature tensor
\begin{equation}
R^{\rho}{}_{\sigma \mu \nu} :=
\frac{\partial}{\partial x^{\mu}} \, \Christ^{\rho}{}_{\nu \sigma} - 
\frac{\partial}{\partial x^{\nu}} \, \Christ^{\rho}{}_{\mu \sigma}
+\Christ^{\rho}{}_{\mu \gamma} \, \Christ^{\gamma}{}_{\nu \sigma} 
-\Christ^{\rho}{}_{\nu \gamma} \, \Christ^{\gamma}{}_{\mu \sigma},
\end{equation}
becomes identically zero.

\section{Under conformal transformation}
\label{appendix:B}
For a Riemannian metric $g$ on a smooth manifold $\mathcal{M}$, and a smooth real-valued function  $\lambda$ on $\mathcal{M}$, 
\begin{equation}
    \tilde{g} = {\rm e}^{2 \lambda} \, g,
\end{equation}
is also a Riemannian metric on $\mathcal{M}$, and we say that $\tilde {g}$ is conformal to $g$.
Under this conformal transformation, we have
\begin{align}
 \tilde {g}_{\mu \nu} &= {\rm e}^{ 2 \lambda } \, g_{\mu \nu}, \\
  \tilde{\Christ}\,{}^{\rho}{}_{\mu \nu} &= \Christ{}^{\rho}{}_{\mu \nu} 
+ \delta_{\mu}^{\rho} \, \partial _{\nu} \lambda 
   +\delta_{\nu}^{\rho} \, \partial _{\mu} \,\lambda -g_{\mu \nu} \nabla ^{\rho} \lambda,
 \label{UnderCT}
\end{align}
where $ \tilde {g}_{\mu \nu}$ is the transformed metric and $\tilde{\Christ}\,{}^{\rho}{}_{\mu \nu} $
is the transformed coefficients of the connection.
It is known that the difference between the Christoffel symbols of two different metrics always form the components of a tensor.



\end{document}